\documentclass[a4paper,11pt]{scrartcl}
\pdfoutput=1
\usepackage[utf8]{inputenc}
\usepackage[T1]{fontenc}
\usepackage{amsmath}
\usepackage[active]{srcltx}
\usepackage{amsthm}
\usepackage{amssymb}
\usepackage{varioref}
\usepackage{makeidx}
\usepackage{longtable}
\usepackage{tikz}
\usetikzlibrary{arrows,matrix,decorations,positioning,backgrounds,fit}

\newtheorem{Def}{Definition}[section]
\newtheorem{Thm}[Def]{Theorem}
\newtheorem{Lemma}[Def]{Lemma}
\newtheorem{Korollar}[Def]{Corollary}
\newtheorem{Bemerkung}[Def]{Remark}

\newtheorem{Vermutung}[Def]{Conjecture}
\newtheorem{Beispiel}[Def]{Example}

\newcommand{\A}{\mathop{\mathbb{A}}}

\newcommand{\C}{\mathbb{C}}

\newcommand{\Z}{\mathbb{Z}}
\newcommand{\N}{\mathbb{N}}
\newcommand{\tr}{\mathrm{tr}}
\newcommand{\Tr}{\mathrm{Tr}}
\newcommand{\GL}{\mathrm{GL}}
\newcommand{\mdeg}{\mathop{\mathrm{mdeg}}}
\newcommand{\M}{\mathrm{M}}

\newcommand{\nh}{\mathrm{N}}

\title{A presentation of the trace algebra of three $3 \times 3$ matrices}
\author{Torsten Hoge}
\date{}
\begin{document}
\maketitle
	\begin{abstract}
		The trace algebra $C_{nd}$ is generated by all traces of 
		products of $d$ generic $n \times n$ matrices. Minimal generating
		sets of $C_{nd}$ and their defining relations are known for $n < 3$ 
		and $n = 3, d=2$.
		This paper states a minimal generating set and their defining relations
		for $n=d=3$. Furthermore the computations yield a description of $C_{33}$ 
		as a free module over	the ring generated by a homogeneous system of 
		parameters.
	\end{abstract}
	\section{Introduction}
	Let $\GL_n:= \GL_n(\C)$ be the general linear group over $\C$ and $\M_n:=\M_n(\C)$
	be the set of $n \times n$ matrices with entries in $\C$. Here $\GL_n$ acts on 
	$\M_n^d$ by simultaneous conjugation, i.e. 
	$g.(A_1,\ldots,A_d) = (gA_1g^{-1},\ldots,gA_dg^{-1})$ for $g \in \GL_n$ and
	$A_1,\ldots,A_d \in \M_n$. This action extends to the coordinate ring $\C[\M_n^d]$ 
	of $\M_n^d$, which is a polynomial ring generated by the projections
	$x_{ij}^{(k)} \colon \M_n^d \rightarrow \C$, 
	$x_{ij}^{(k)}(A_1,\ldots,A_d) = (A_k)_{ij}$. The action is given by 
	$(g.x_{ij}^{(k)})(A) = x_{ij}^{(k)}(g^{-1}.A)$. The invariant ring
	\begin{displaymath}
		\C[\M_n^d]^{\GL_n} = \{ f \in \C[\M_n^d] \mid g.f = f \text{ for all } g \in 
		\GL_n \}
	\end{displaymath}
	is the coordinate ring of the algebraic quotient of the above action of $\GL_n$
	on $\M_n^d$. We want	to find a minimal presentation of this algebra, i.e. a 
	minimal homogeneous	generating set and a minimal set of defining relations 
	between these generators. 

	The next theorem shows that the invariant ring is generated by the traces
	of products of generic matrices.
	\begin{Thm}[First fundamental theorem \cite{procesi}]
	\begin{displaymath}
		C_{nd}:= \C[\M_n^d]^{\GL_n} = \C[\tr(X_{i_1} \cdots X_{i_k}) \mid i_1,\ldots,
			i_k \in \{1,\ldots,d\}, k \in \N],
	\end{displaymath}
	where $X_k = (x_{ij}^{(k)})$ is the $k$-th generic matrix formed by the 
	projections. 
	\end{Thm}
	From this description it is obvious that the invariant ring
	is a (multi-)graded algebra. Here the degree is given by 
	$\deg \tr(X_{i_1} \cdots X_{i_k}) = k$ which	is compatible with the standard 
	grading of $\C[\M_n^d]$. The multigrading counts	how often the generic 
	matrices occur, for example $\mdeg \tr(X_1X_3^2X_1X_3X_2) = (2,1,3)$. 

	This algebra is finitely generated (see \cite{mumford:git}§2) and there is a 
	degree bound 	$\nh(n)$ for the generators which is sharp for large $d$ 
	(\cite{procesi}).
	The bound is given by
	the nilpotency class of certain (non-unitary) algebras in the Nagata-Higman 
	theorem	\cite{procesi}. The values of this function are only known for $n<5$ 
	($\nh(1)=1$, $\nh(2)=3$, $\nh(3)=6$, $\nh(4)=10$). 
	The upper bound 	$\nh(n) \le n^2$ for all $n$ was 
	given by Razmislov (see \cite{drensky:pir} 6.2 and \cite{formanek:nagatahigman} 
	for the details about $\nh(n)$). 
	This yields a simple algorithm to compute a minimal homogenous generating set. 
	Namely we write down the traces mentioned in the first fundamental theorem 
	up to degree $\nh(n)$ and check if they are generated by the other ones.
	There is a more sophisticated possibility given in \cite{drensky:4x4gen} 
	and \cite{drensky:3x3def} which uses a $\GL_d$-action on $C_{nd}$. 	The 
	minimal generating set for $C_{33}$ given in \ref{generators} was computed 
	in this way. The generators are grouped by the $\GL_d$-action on $C_{nd}$. 
  The first generator	in each group is the highest weight vector with respect 
	to this action.	We refer to these generators by $t_k$ where $k$ is given in 
	the list \ref{generators}.
	We fix this	minimal generating set of $C_{33}$ for the rest of the article.

	Given such a minimal homogeneous generating set there is a minimal graded
	free resolution 
  \begin{displaymath}
  0 \rightarrow \bigoplus_j R[-j]^{\beta_{kj}} \rightarrow \ldots \rightarrow 
	\bigoplus_j R[-j]^{\beta_{1j}} \rightarrow R \rightarrow C_{nd} \rightarrow 0
  \end{displaymath}
	of $C_{nd}$ (\cite{bruns:cmr} chapter 1.5). Here $R = \C[T_1,\ldots,T_k]$ 
	where the $T_i$ correspond to the elements of the minimal generating set.
	And $R[-j]$ denotes the shift of the grading of $R$ by $j$ so that $1 \in R[-j]$ 
	has degree $j$.
	Furthermore the $\beta_{ij}$ do not depend on the choice of 
	the minimal homogeneous generating set. The map 
	$\bigoplus_j R[-j]^{\beta_{1j}} \rightarrow R$ decodes the minimal generating 
	set of the ideal of relations among the generators. 
	Similar to the Nagata-Higman bound we define
  \begin{displaymath}
   \nh^i(n,d) := max_j \{\beta_{ij} \not=0\}.
  \end{displaymath}

	There a two simplifications mentioned already in \cite{drensky:3x3def}. First 
	\begin{displaymath}
		C_{nd} = \C[\tr(X_1),\ldots,\tr(X_d)] \otimes_\C \C[Q]
	\end{displaymath}
	where $Q$ is generated by generic matrices with trace $0$. One gets this
	description by $x_i := X_i -n\,\tr(X_i)E$ where $E$ is the $n \times n$ 
	identity matrix.
	We will denote these generic traceless matrices
	with small letters. By the description above the relations of $C_{nd}$ are
	given by the relations of $\C[Q]$.
  The second simplification is that we can assume that $x_1$ is a diagonal 
	generic matrix, because the diagonalizable matrices form a dense subset of 
	$M_n$. This is essential	for the computations, because we get rid of 
	$d + n(n-1)$ variables.
  In case $d=3$ we have only $3$ generic matrices.
	We denote these by $X,Y,Z$ and the traceless ones by $x,y,z$. 

	We have seen that $C_{nd}$ is a graded algebra. But it also has the property
	of beeing a Cohen-Macaulay ring (\cite{drensky:computing}). This means there 
	is a sequence of homogeneous sucessive non-zero divisors of length 
	$\dim C_{nd}$ called a maximal homogeneous regular sequence.	
	For a positive graded Cohen-Macaulay ring 
	maximal homogeneous regular sequences coincide with homogeneous systems of 
	parameters (see \cite{springer}). This will be used to compare the Hilbert series of 
	$C_{nd}$ with the Hilbert series given by the generators and relations. The
	Hilbert series of $C_{nd}$ was computed in \cite{berele:hilbert} for small cases, 
	in particular for $n=3=d$.

	If we fix a minimal homogeneous generating set $E = \{t_1,\ldots,t_k\}$ of 
	$C_{nd}$ we write $\C_E = \C[T_1,\ldots,T_k]$ for the corresponding polynomial
	ring.	The $T_i$ correspond to the generators and we have a canonical map onto
	$C_{nd}$ given by	$T_i \mapsto t_i$. 

	\section{Algorithms}\label{algorithms}
	Given a minimal homogeneous generating set $\{t_1,\ldots,t_k\}$ of $C_{nd}$ we 
	get a morphism of affine varieties $\varphi \colon \M_n^d \rightarrow \A^k$ with
	$\varphi(A) = (t_1(A),\ldots,t_k(A))$. The comorphism
	$\varphi^\ast \colon \C[T_1,\ldots,T_k] \rightarrow \C[\M_n^d]$ sends
	$T_i \mapsto T_i \circ \varphi = t_i$. The kernel of this comorphism is the ideal
	of relations of the minimal generating set. If we fix the degrees of the $T_i$ 
	by $\deg T_i := \deg t_i$ the comorphism is a graded algebra homomorphism.
	So we can compute the kernel degree by degree. In each degree we have to solve
	one linear system (which is given by the coefficients of the images of the
	monomials of that degree in $\C[T_1,\ldots,T_k]$). The problem here is that 
	these linear systems grow quite fast and furthermore it is not obvious how 
	many degrees we have to consider. We will see how one could obtain a degree
	bound by a theorem of Harm Derksen in \ref{upper bound}, which uses the 
	Cohen-Macaulay property	of $C_{nd}$. 

	The second algorithm is a consequence of the second fundamental theorem. This 
	theorem describes the relations of $C_{nd}$ in terms of the generators which 
	occur in the first fundamental theorem. The formal trace algebra $C_\infty$
	is the polynomial ring generated by formal traces, i.e. it is generated by 
	formal traces $\Tr(w)$ where $w$ is a word in $X_1,X_2,\ldots$. But we have
	to consider the trivials relation between traces. Two formal traces $\Tr(w)$ 
	and $\Tr(w')$ are equal	if and only if $w$ is a cyclic permutation of $w'$.
	We can also consider the formal trace algebra generated by $d$ letters and
	denote this algebra by $C_{\infty,d}$. This algebra is graded in the same
	way as the grading defined on $C_{nd}$.
	We have a canonical map
	$\pi\colon C_\infty \rightarrow C_n$ which replaces $\Tr$ by $\tr$ and
	the letters by generic $n \times n$ matrices. Here 
	\begin{displaymath}
		C_{n}:=  \C[\tr(X_{i_1} \cdots X_{i_k}) \mid i_1,\ldots,
			i_k,k \in \N]
	\end{displaymath}
	is given as in the first fundamental theorem but we allow arbitrary many
	generic matrices.

	\begin{Def}
		A \emph{trace identity} for $n \times n$ matrices is an element 
		$f \in C_\infty$ such that $\pi(f) = 0$.
	\end{Def}

	\begin{Def} [\cite{procesi}]The fundamental trace identity for $n \times n$ matrices is given by
		\begin{displaymath}
			F(X_1,\ldots,X_{n+1}) := \sum_{\sigma \in S_{n+1}} \Tr_\sigma(X_1,\ldots,X_{n+1}).
		\end{displaymath}
		Here $S_{n+1}$ is the permutation group of ${n+1}$ elements.
	\end{Def}
	 Example: $\Tr_{(12)(3)}(X,Y,Z) = \Tr(XY)\Tr(Z)$. 
	\begin{Thm}[Second fundamental theorem \cite{procesi}]
		Every trace identity for $n\times n$ matrices is included in the ideal 
		generated by the
		\begin{center}
			$F(M_1,\ldots,M_{n+1})$,
		\end{center}
		where the $M_i$ are (non-constant) monomials in the $X_i$.
	\end{Thm}

	The advantage of the second fundamental theorem is that a relation can be
	given as a tuple of monomials. This is a very efficient way to describe the 
	generators of the ideal. The relations of $C_{nd}$ are quite complicated. The
	goal in this section is to describe the minimal generating set of the 
	relations of $C_{33}$ by some of these monomials and some extra data.
	
	We define $I_{\infty,d}$ to be the ideal in $C_{\infty,d}$ which is
	generated by the $F(M_1,\ldots,M_{n+1})$ where only the letters 
	$X_1,\ldots,X_d$ occur in the monomials. The ideal $I_{\infty,d}$ is generated 
	by formal traces of arbitrary degree.
	On the other hand	the relations according to the minimal homogeneous 
	generating set of $C_{nd}$ are also part of this ideal. The main problem is 
	to rewrite the formal traces in terms of the minimal generating set. 
	Here $C_E = \C[T_1,\ldots,T_k]$ can be seen as a subalgebra of $C_{\infty,d}$.
	So rewriting is given by an evaluation map 
	$R \colon C_{\infty,d} \rightarrow \C[T_1,\ldots,T_k]$. 
	Here we only consider	evaluation maps which are constant on 
	$\C[T_1,\ldots,T_k]$. Furthermore the diagram
	\begin{center}
		\begin{tikzpicture}[description/.style={fill=white,inner sep=2pt}]
	  \matrix (m) [matrix of math nodes, row sep=3em, column sep=2.5em, text height=1.5ex, text depth=0.25ex]
	  {	C_{\infty,d} & C_{nd} \\ 
			C_E& C_{nd} \\};
	  \path[->, font=\scriptsize]
		  (m-1-1) edge node[auto] {$\pi$} (m-1-2)
		  (m-1-1) edge node[auto] {$R$} (m-2-1)
		  (m-2-1) edge node[auto] {$\pi|_{C_E}$} (m-2-2)
		  (m-1-2) edge node[auto] {$id$} (m-2-2);
		\end{tikzpicture}
	\end{center}	
	should commute and $R$ should be a (multi-)graded map. These are quite
	natural assumptions. From the diagram	we see that $R(x)$ and $x$ may differ 
	up to an element of $\ker \pi$. 

	The	next lemma shows that we only have to fix an evaluation map and a set of 
	tuples $M = (M_1,\ldots,M_{n+1})$ as in the second fundamental theorem.
	\begin{Lemma}
		Let $R\colon C_{\infty,d} \rightarrow C_E$ be an graded evaluation map 
		which is a projection on $C_E$ and commutes with $\pi$. Then there exist 
		a set $S$ of tuples of monomials such that the set $S_R:=\{R(F(M))\mid M \in S\}$ 
		is a minimal generating set of $\ker(\pi|_{C_E})$.
		\begin{proof}
			Let $r \in \ker(\pi|_{C_E})$. Then $r \in I_{\infty,d}$. From the second fundamental 
			theorem follows that there are 	$c_M \in C_{\infty,d}$ with	
			$r= \sum_M c_M F(M)$.
		  Then $r= R(r) = \sum_M R(c_M) R(F(M))$. So we get a finite generating set
			of the relations by choosing all $R(F(M))$ up to degree $N^1(n,d)$. 
			Since $C_{nd}$ is a positive graded ring and the $R(F(M))$ are homogeneous 
			we can choose a minimal generating set which consists of some of the 
			$R(F(M))$.
		\end{proof}
	\end{Lemma}
	
	There is one problem with this setup. Take an element $F(M)$ with 
	$R(F(M)) \not= 0$. Because $R(F(M)) \in I$ we can define another evaluation 
	$R'$ by just changing the evaluation of the formal traces of maximum degree 
	in $F(M)$ in the following way. 
	Let $R'(c):= R(c) - \frac{1}{n!} R(F(M))$ for these formal traces. 
	Then $R'(F(M)) = 0$ and this is not a part of a minimal 
	generating set. That means the tuples $M$ depend on the evaluation map.

	So we have to define such an evaluation map. One possibility would be 
	to compute the evaluation for every element by solving a system of linear 
	equations. But then we could also use the first method and furthermore we 
	would only fix it up to a given degree.	The second method depends on the 
	next definition. 

	\begin{Def}\label{trace reduction}
		A \emph{trace reduction for $n \times n$ matrices} is a multi-homogeneous 
		trace identity of $n \times n$ matrices
		\begin{equation}
		 \Tr(X_1X_2\cdots X_{\nh(n)} X_{\nh(n)+1}) = \sum_w \lambda_w \prod_{w} \Tr(w),
		\end{equation}
		where $w$ are words in the letters $\{X_1,\ldots,X_{\nh(n)+1}\}$ and 
		$\lambda_w \in \C$.
		Trace identity means that
		\begin{equation}
		 \Tr(X_1X_2\cdots X_{\nh(n)} X_{\nh(n)+1}) - 
		 \sum_w \lambda_w \prod_{w} \Tr(w) \in \ker \pi.
		\end{equation}
	\end{Def}

	It is obvious from the definition	how one could reduce a trace of 
	degree $\nh(n)+1$ in terms of traces of lower degrees. 
	Since all the variables of $C_\infty$ correspond to	such traces we can 
	reduce them in a unique way by the following algorithm. Here we write
	$\Tr(Y_{i_1} \cdots Y_{i_k})$ instead of $\Tr(X_{i_1} \cdots X_{i_k})$ for 
	the elements we want to reduce to distinguish them from the trace reduction.
	\begin{enumerate}
		\item A generator of $C_\infty$ corresponds to a formal trace 
			$Tr(Y_{i_1} \cdots Y_{i_k})$. 
		\item If we permute these $i_j$ by a cyclic permutation the image under
			$\pi$ is the same. All the elements in that orbit should be reduced
			to the same element. So choose the one where the defining tuple 
			$(i_1,\ldots,i_k)$ is minimal for the lexicographical order.
		\item Replace this formal trace with the right hand side of the trace 
			reduction with $X_1:= Y_{i_1},\ldots,X_{\nh(n)}:=Y_{i_{\nh(n)}}$ and 
			$X_{\nh(n)+1}:= Y_{i_{\nh(n)+1}}\cdots Y_{i_k}$.
		\item Repeat this until all given traces have degree $\le \nh(n)$.
		\item Reduce traces in terms of the choosen minimal generating set 
					(unique for $n=3$).
	\end{enumerate}

	The next question is how one gets such a trace reduction. 
  This can be done by the second fundamental theorem. We know
	that $\Tr(X_1\cdots X_{\nh(n)+1})$ only occurs in $F(M)$ if and only if 
	each letter $X_1,\ldots,X_{\nh(n)+1}$ in the monomials occur exactly once. 
	When one inserts such	a monomial into the fundamental trace identity 
	it is ovious that only the	traces of full degree which
	are not $\Tr(X_1X_2\cdots X_{\nh(n)} X_{\nh(n)+1})$ have to be eliminated. 
	Let $m$ be the number of tuples of monomials which fulfill the above 
	assumptions.	We get a $\nh(n)! \times m$-matrix $A$ with 
	entries in $\Z$ (evidently all entries are $0$ or $1$) where the
	rows are indexed by the formal traces of length $\nh(n)+1$ and the columns are 
	indexed by these tuples of monomials. If we order these traces
	such that $\Tr(X_1\cdots X_{\nh(n)+1})$ is the first one, we only have
	to solve the linear equation $Ax=e_1$. Then the entries of $x$ are the 
	coefficients of the corresponding trace reduction.
			
	\begin{Beispiel} 
		We find the trace reduction for $2 \times 2$-matrices. Here $\nh(2)=3$ so
		we need to find a reduction for $\Tr(X_1X_2X_3X_4)$. If we evaluate 
		\begin{eqnarray*}
		F(X_1X_2,X_3,X_4) &=& \Tr(X_1X_2X_3X_4) + \Tr(X_1X_2X_4X_3) 
		- \Tr(X_1X_2)\Tr(X_3X_4)\\ 
		&-& \Tr(X_1X_2X_3)\Tr(X_4) - \Tr(X_1X_2X_4) \Tr(X_3)	\\
		&+& \Tr(X_1X_2)\Tr(X_3)\Tr(X_4)
		\end{eqnarray*}
		we see that only the first two terms on the right hand side are relevant.
		Because the notation is quite clumsy we will denote $\Tr(X_1X_2)\Tr(X_3X_4)$ by
		$[12][34]$, which is not ambigious as long as we restrict to less than
		$10$ matrices. We will further write $F([12,3,4])$ for the left hand side.
		Now we get some equations:
	  \begin{eqnarray*}
	   F([12,3,4]) &\hat=& [1234] + [1243] \\
	   F([41,2,3]) &\hat=& [1234] + [1324] \\
	   F([24,1,3]) &\hat=& [1324] + [1243].\\
	  \end{eqnarray*}		
		Here we write $\hat=$ because we only consider the traces of maximal length.
		This gives rise to the following linear system
		\begin{displaymath}
		\begin{pmatrix}
			1&1&0 \\
			1&0&1 \\
			0&1&1 \\
		\end{pmatrix}
		x = 		
		\begin{pmatrix}
			1  \\
			0 \\
			0 \\
		\end{pmatrix}
		\end{displaymath}
		with solution $x = \frac{1}{2}\begin{pmatrix}
			1  \\
			1 \\
			-1 \\
		\end{pmatrix}
		$. This means the trace reduction is given by
   \begin{eqnarray*}
    [1234]&=&  - \frac{1}{2} ( [243][1] - [123][4] -[124][3] - [134][2] 
							 + [12][3][4]  \\
    &&  \quad + [14][2][3]- [24][1][3] - [12][34] - [14][23] + [24][13]). \\
   \end{eqnarray*}
		We have already seen that we may consider traceless matrices. So we may
		assume that $X_1,X_2$ and $X_3$ correspond to traceless matrices and denote
		them by $x_1,x_2,x_3$. In the reduction we will insert products of matrices
		into $X_4$, so $X_4$ we cannot assume $X_4$ to be traceless. In this case 
		we get a more compact	trace reduction
  	\begin{eqnarray} \label{red22}
	    \Tr(x_1x_2x_3X_4) 	&=& \frac{1}{2} (\Tr(x_1x_2x_3)Tr(X_4) + \Tr(x_1x_2)\Tr(x_3X_4) \\
			& & + \Tr(x_1X_4)\Tr(x_2x_3) - \Tr(x_2X_4)\Tr(x_1x_3)). \notag
  	\end{eqnarray}
		Using traceless matrices has the advantage that we always reduce the degree 
		by at least $2$ (at least if we are not in degree $\nh(n)+1$). We could enlarge 
		the linear systems with terms like $F([23,1,4])$ such that the matrix
		has a non-trivial kernel. But there is only one trace reduction if we
		restrict ourselfes to terms of a minimal generating set of $C_{2d}$, i.e.
		$\tr(X_1),\ldots,\tr(X_d),\tr(X_iX_j) i\le j, \tr(X_iX_jX_k)$ for $i <j<k$. 
		With these assumptions there is only one trace reduction, because there are 
		no non-trivial relations of degree $4$ in $C_{2d}$ see
		(\cite{drensky:pir} Theorem 5.3.8) and an 
		element	of the kernel of the above matrix would give rise to such an 
		non-trivial relation.
	\end{Beispiel}

	\begin{Bemerkung}
		The example shows how one can get such a trace reduction for $n\times n$
		matrices. Write down all coefficients of the cycles in $F(M)$ for every 
		suitable $M$, i.e. all numbers from $1$ to $\nh(n)+1$ occur once in $M$.
		Pick a solution of the corresponding linear system. This always works
		since these formal traces of degree $\nh(n)+1$ only occur for such $M$ isolated
		and by the definition of $\nh(n)$ there has to be such a relation.
		Unfortunately such a trace reduction is not unique for $n \ge 3$.
	\end{Bemerkung}	

  The set of trace reductions
	depends on the kernel of a linear map, i.e. we can describe them as an
	affine variety. The reduction of $R(F(M))$ only depends on the 
	coefficients of this affine variety. That such a set of $R(F(M))$ is 
	a minimal generating set is equivalent to 
	$\dim (R(F(M))_k = \dim (\ker(\pi|_{C_E})_k)$ for
	a finite number of degrees $k$ (It is enough to check the degrees of the 
	defining relations). If the set of $R(F(M))$ are not a generating set then 
	a $\dim (\ker(\pi|_{C_E})_k)$	minor of the corresponding set of linear equations 
	given by the $R(F(M))$ vanishes. Since we only have to check a finite number 
	of degrees, the product	of the corresponding minors gives us the 
	equation, when such a trace reduction	eliminates the minimal generating set 
	given by the fixed set of tuples of monomials.

	\begin{Lemma}
		Let $R$ be an evaluation map which is given by a trace reduction. 
		Further let S be a (minimal) set of tuples of monomials such that 
		$S_R=\{R(F(M))\mid M \in S\}$ is a minimal generating set of 
		$\ker (\pi|_{C_E})$. 
		Then the set 
		\begin{displaymath}
			\{R' \mid R' \text{ given by trace reduction with } (S_{R'}) = \ker(\pi|_{C_E})\}
		\end{displaymath}
		is generic.
		\begin{proof}
			The set of $R'$ with $(S_{R'}) \not= \ker(\pi|_{C_E})$ is closed by 
			the observation above. So the complement is a dense subset.
		\end{proof}
	\end{Lemma}		

	\begin{Bemerkung}
		The lemma tells us that it is unlikely to choose the wrong trace 
		reduction. Therefore the main theorem \ref{main theorem} only states the
		tuples of monomials. One can find the corresponding 
		trace reduction in \cite{doktorarbeit} (page 95-103).
	\end{Bemerkung}

	\begin{Vermutung}
		If we choose $R$ as in \ref{trace reduction} then the choosen monomials
		always form a minimal generating set.
	\end{Vermutung}

	\section{An upper bound for the relations}\label{upper bound}
	From \cite{derksen:upper_bounds} follows that the degree of the defining 
	relations of $C_{nd}$	is bounded by
	\begin{displaymath}
		d_1 + d_2 + \ldots + d_{\dim C_{nd}+1} + a(C_{nd}).
	\end{displaymath}
	Here the $d_i$ are the degrees of the elements of the minimal generating
	set ordered descending and $a(C_{nd})$ is the degree of the Hilbert series
	of $C_{nd}$ (degree of nominator - degree of denominator).
	Here $a(C_{nd}) \le -\dim C_{nd}$ holds due to 
	\cite{derksen:upper_bounds} via \cite{Knop:kanonisch}. 
	Since $\dim C_{nd} = (d-1) n^2 +1$ (\cite{lopatin}) for $d \ge 2$ and 
	the generators are bounded by $\nh(n) \le n^2$ one gets an upper bound for the 
	defining relations
	\begin{displaymath}
		\nh^1 (n,d) \le ((d-1)n^2+2)n^2 - ((d-1)n^2+1)= (d-1)n^4+(3-d)n^2 -1.
	\end{displaymath}
	For $C_{33}$ follows $\nh^1(3,3) \le 161$ which is quite a bad bound.
	With some more concret values for $C_{33}$, i.e. $a(C_{33}) = -27$,
	$(d_i) = (6^{10},5^9,4^9,3^{11},2^6,1^3)$ one gets
	\begin{displaymath}
		\nh^1(3,3) \le 82.
	\end{displaymath}

	This is also not a very sharp bound. The following lemma allows us to 
	examine the whole setup modulo a homogeneous system of parameters.

	\begin{Lemma}
		let $R=\C[X_1,\ldots,X_k]$ be the polynomial ring in $k$ variables and
		$\deg X_i >0$ for all $i$. Further let $I$ be a homogeneous ideal of $R$ and
		$\{f_1,\ldots,f_n\}$ be a homogeneous minimal generating set of $I$ with 
		$\deg f_i >0$ for all $i$. If $g \in R$
		is a homogeneous non-zerodivisor of $R/I$ with $\deg g > 0$ then 
		$\{f_1,\ldots,f_n,g\}$ is a minimal	homogeneous generating set of the 
		ideal $(I,g)$.
		\begin{proof}
			First $g \notin I$ because $g$ is $R/I$-regular. It is enough
			to show that no $f_i$ can be expressed by $g$ and the other $f_j$. Assume
			\begin{displaymath}
				f_i = \sum_{j\not=i} \beta_j f_j + \beta g
			\end{displaymath}
			where $\beta,\beta_j \in R$. So $\beta g \in I$ and therefore 
			$\beta \in I$ because $g$ is a non-zerodivisor of $R/I$. So there exist
			$\alpha_k\in R$ with $\beta = \sum_{\deg(f_k) < \deg(f_i)} \alpha_k f_k$ 
			because $g$ has	positive degree. But now
			\begin{displaymath}
				f_i = \sum_{j\not=i} \beta_j f_j + 
				\sum_{\deg(f_k) < \deg(f_i)} g \alpha_k f_k
			\end{displaymath}
			can be expressed by the other $f_i$. This is a contradiction to the 
			minimality of $\{f_1,\ldots,f_n\}$.
		\end{proof}
	\end{Lemma}

	Using this lemma for the homogeneous parametersystem $H$ given in \ref{hompar}
	we get
	\begin{displaymath}
		\nh^1(3,3) \le 6 -27 + 48 = 27.
	\end{displaymath}
	Here $a(R/(I,H))=48-27$. 
\pagebreak
\begin{Lemma} \label{generators}The following elements form a minimal generating
set of $C_{33}$. \\
\begin{tabular}{ll}
\underline{$W_3(1)$:} & \underline{$W_3(2^2)$:}\\
$(a)$ $\tr(X)$ & $(18)$ $\tr(x^2y^2)-\tr(xyxy)$\\
$(b)$ $\tr(Y)$ & $(19)$ $\tr(x^2z^2)-\tr(xzxz)$\\
$(c)$ $\tr(Z)$ & $(20)$ $\tr(y^2z^2)-\tr(yzyz)$\\
							& $(21)$ $\tr(x^2yz) + \tr(x^2zy)- 2 \;\tr(xyxz)$\\
\underline{$W_3(2)$:}& $(22)$ $\tr(y^2xz) + \tr(y^2zx)- 2 \;\tr(yxyz)$\\
$(1)$ $\tr(x^2)$ & $(23)$ $\tr(z^2xy)+\tr(z^2yx)-2 \;\tr(zxzy)$\\
$(2)$ $\tr(y^2)$ & \\
$(3)$ $\tr(z^2)$ & \underline{$W_3(2,1^2)$:}\\
$(4)$ $\tr(xy)$ & $(24)$ $\tr(x^2yz)-\tr(x^2zy)$\\
$(5)$ $\tr(xz)$ & $(25)$ $\tr(y^2xz)-\tr(y^2zx)$\\
$(6)$ $\tr(yz)$ & $(26)$ $\tr(z^2xy)-\tr(z^2yx)$\\
&\\
\underline{$W_3(3)$:}& \underline{$W_3(3,1^2)$:}\\
$(7)$ $\tr(x^3)$& $(27)$ $\tr(x^2yxz)-\tr(x^2zxy)$\\
$(8)$ $\tr(y^3)$& $(28)$ $\tr(y^2xyz)-\tr(y^2zyx)$\\
$(9)$ $\tr(z^3)$& $(29)$ $\tr(z^2yzx)-\tr(z^2xzy)$\\
$(10)$ $\tr(x^2y)$& $(30)$ $\tr(yxyxz)+\tr(x^2y^2z)-\tr(xyzxy)-\tr(x^2zy^2)$\\
$(11)$ $\tr(x^2z)$& $(31)$ $\tr(zxzxy)+\tr(x^2z^2y)-\tr(xzyxz)-\tr(x^2yz^2)$\\
$(12)$ $\tr(y^2x)$& $(32)$ $\tr(yzyzx)+\tr(z^2y^2x)-\tr(zyxzy)-\tr(z^2xy^2)$\\
$(13)$ $\tr(y^2z)$&\\
$(14)$ $\tr(z^2x)$&\\
$(15)$ $\tr(z^2y)$&\underline{$W_3(1^3)$:}\\
$(16)$ $\tr(xyz)+\tr(xzy)$&$(17)$ $\tr(xyz)-\tr(xzy)$\\
&\\
\underline{$W_3(2^2,1)$:}&\\
\multicolumn{2}{l}{$(33)$ $\tr(x^2y^2z)+\tr(x^2zy^2)+\tr(xyxyz)+\tr(xyxzy)-2\tr(x^2yzy)-2\tr(xy^2xz)$}\\
\multicolumn{2}{l}{$(34)$ $\tr(x^2z^2y)+\tr(x^2yz^2)+\tr(xzxzy)+\tr(xzxyz)-2\tr(x^2zyz)-2\tr(xz^2xy)$}\\
\multicolumn{2}{l}{$(35)$ $\tr(y^2z^2x)+\tr(y^2xz^2)+\tr(yzyzx)+\tr(yzyxz)-2\tr(y^2zxz)-2\tr(yz^2yx)$}\\
&\\
\underline{$W_3(3^2)$:}&\\
\multicolumn{2}{l}{$(36)$ $\tr(x^2y^2xy)-\tr(y^2x^2yx)$}\\
\multicolumn{2}{l}{$(37)$ $\tr(x^2z^2xz)-\tr(z^2x^2zx)$}\\
\multicolumn{2}{l}{$(38)$ $\tr(y^2z^2yz)-\tr(z^2y^2zy)$}\\
\multicolumn{2}{l}{$(39)$ \small{$\tr(x^2yxyz)+\tr(x^2yxzy)+\tr(x^2zxy^2)-\tr(x^2y^2xz)-\tr(x^2yzxy)-\tr(x^2zyxy)$}}\\
\multicolumn{2}{l}{$(40)$ \small{$\tr(y^2xyxz)+\tr(y^2xyzx)+\tr(y^2zyx^2)-\tr(y^2x^2yz)-\tr(y^2xzyx)-\tr(y^2zxyx)$}}\\
\multicolumn{2}{l}{$(41)$ \small{$\tr(z^2yzyx)+\tr(z^2yzxy)+\tr(z^2xzy^2)-\tr(z^2y^2zx)-\tr(z^2yxzy)-\tr(z^2xyzy)$}}\\
\multicolumn{2}{l}{$(42)$ \small{$\tr(x^2zxzy)+\tr(x^2zxyz)+\tr(x^2yxz^2)-\tr(x^2z^2xy)-\tr(x^2zyxz)-\tr(x^2yzxz)$}}\\
\multicolumn{2}{l}{$(43)$ \small{$\tr(y^2zyzx)+\tr(y^2zyxz)+\tr(y^2xyz^2)-\tr(y^2z^2yx)-\tr(y^2zxyz)-\tr(y^2xzyz)$}}\\
\multicolumn{2}{l}{$(44)$ \small{$\tr(z^2xzxy)+\tr(z^2xzyx)+\tr(z^2yzx^2)-\tr(z^2x^2zy)-\tr(z^2xyzx)-\tr(z^2yxzx)$}}\\
\multicolumn{2}{l}{$(45)$ \small{$\tr(x^2y^2z^2)+\tr(x^2zyzy)+\tr(xyxz^2y)+\tr(xyzxyz)+\tr(xzxzy^2)$}}\\
\multicolumn{2}{l}{\quad \, \small{$-\tr(x^2yzyz)-\tr(x^2z^2y^2)-\tr(xyxyz^2)-\tr(xy^2zxz)-\tr(xzyxzy)$}}\\
\end{tabular}
\end{Lemma}

\section{The relations}
\begin{Thm} \label{main theorem}
	The following tuples represent the minimal set of relations between the 
	generators given in \ref{generators}
\begin{center}
\begin{tabular}{cccc}
  \bf{degree $7$}			&  										&										 												\\  
  \bf{deg $(3,2,2)$}	& \bf{deg $(3,3,3)$}	&	\bf{deg $(4,3,3)$} 	&	\bf{deg $(4,4,3)$}	\\  
	$(111,22,3,3)$ 			&	$(1132,223,1,3)$		&	$(1111,222,33,3)$		& $(1212,2121,33,3)$	\\
								 			&	$(1332,223,1,1)$		&	$(1111,222,33,3)$		& $(1212,2112,33,3)$	\\
											& $(3332,221,1,1)$		&	$(1112,122,33,3)$		& $(1212,2211,33,3)$	\\  
  \bf{degree $8$}		 	& $(3331,221,2,1)$		&	$(1122,112,33,3)$		& $(1212,1221,33,3)$	\\  
  \bf{deg $(4,3,1)$}	&	$(1112,223,3,3)$		&	$(1322,112,13,3)$		& $(1212,1212,33,3)$	\\
	$(1111,22,2,3)$		 	& $(2223,331,1,1)$		&	$(22213,11,13,3)$		& $(1212,1122,33,3)$	\\  
  \bf{deg $(4,2,2)$}	& $(2213,331,2,1)$ 		&	$(22123,11,13,3)$		& $(1122,1122,33,3)$	\\
	$(1111,22,3,3)$			& $(111,222,33,3)$		&	$(21223,11,13,3)$		& $(1122,2211,33,3)$	\\
	$(1112,21,3,3)$			& $(112,223,33,1)$ 		& $(12223,11,13,3)$		& $(1122,2121,33,3)$	\\
	$(1122,11,3,3)$			& $(122,231,33,1)$ 		& $(22231,11,13,3)$		& 										\\		
  \bf{deg $(3,3,2)$}	& $(122,231,33,1)$ 		&	$(23221,11,13,3)$		& \bf{degree $12$}		\\
	$(1112,22,3,3)$			&	$(123,123,12,3)$		&	$(32221,11,13,3)$		& \bf{deg $(6,6,0)$}	\\	
	$(1122,21,3,3)$ 		& $(132,123,12,3)$		&	$(22321,11,13,3)$		& $(2112,121,122,12)$	\\
	$(1122,23,1,3)$			& $(132,132,12,3)$		&	$(23123,11,12,3)$		& \bf{deg $(6,5,1)$}	\\
	$(1322,23,1,1)$			&  										&											& $(112212,1123,1,2)$	\\		
	$(3312,11,2,2)$			& \bf{degree $10$}		& \bf{degree $11$}		& \bf{deg $(6,4,2)$}	\\
											&	\bf{deg $(6,2,2)$}	&	\bf{deg $(6,4,1)$}	& $(112122,1133,1,2)$	\\
	\bf{degree $9$}			&	$(121132,11,1,3)$		&	$(1111,1122,22,3)$	& $(112312,1123,1,2)$	\\
	\bf{deg $(5,2,2)$}	&	\bf{deg $(5,4,1)$}	&	\bf{deg $(6,3,2)$}	& \bf{deg $(6,3,3)$}	\\
	$(1111,212,3,3)$		& $(121122,21,1,3)$		& $(1113,1122,21,3)$	& $(112212,1133,1,3)$	\\
	$(1112,112,3,3)$		& $(122122,11,1,3)$		& \bf{deg $(5,5,1)$}	& $(112312,1132,1,3)$	\\
	\bf{deg $(4,4,1)$}	& \bf{deg $(5,3,2)$}	& $(1112,1122,22,3)$	& \bf{deg $(5,2,2)$}	\\
	$(1112,122,2,3)$		&	$(122121,11,3,3)$		&	$(2112,1122,12,3)$	& $(1122,112,122,33)$	\\
	$(1132,122,2,1)$		& $(122111,21,3,3)$		& \bf{deg $(5,4,2)$}	& $(123,132,112,221)$	\\
	\bf{deg $(5,3,1)$}	& $(122311,21,1,3)$		& $(13112,1222,1,3)$	& \bf{deg $(5,4,3)$}	\\
	$(1112,112,2,3)$		&$(122131,21,1,3)$		& $(1133,1122,22,1)$	& $(33131,121,212,2)$	\\
	\bf{deg $(4,3,2)$}	&$(121231,21,1,3)$		& $(1133,1212,22,1)$	& $(1231,132,321,21)$	\\
	$(1111,222,3,3)$		&\bf{deg $(4,4,2)$}		& $(1311,121,222,3)$	& $(123,132,321,112)$	\\
	$(1112,122,3,3)$		&$(123212,21,1,3)$		& $(111,112,222,33)$	& \bf{deg $(4,4,4)$}	\\
	$(1122,112,3,3)$		&$(123212,11,2,3)$		& \bf{deg $(5,3,3)$}	& $(11232,123,123,3)$	\\
	$(1222,113,1,3)$		&$(323212,11,2,1)$		& $(11223,112,13,3)$	& $(13233,11223,1,2)$	\\
	$(3222,111,1,3)$		&$(12211,122,3,3)$		& $(11223,121,13,3)$	& $(13232,1123,13,2)$	\\
	$(3222,111,1,3)$		&$(12121,122,3,3)$		& $(11223,221,13,3)$	& $(1323,123,123,12)$	\\
	$(3322,111,1,2)$		&$(11221,122,3,3)$		& $(11232,112,13,3)$	&											\\
	$(3322,121,1,1)$		&$(11122,122,3,3)$		& $(11232,121,13,1)$	& 										\\
											&$(22111,122,3,3)$		&											&											\\
											&$(22113,122,1,3)$		&											&											\\	
											&$(22311,122,1,3)$		&											&											\\

\end{tabular}
\end{center}
\end{Thm}
\begin{Bemerkung}
	Observe that only the relations are given for the multidegrees which are
	partitions. One gets the other ones by permuting the letters in the formal
	traces and reduce them afterwards in terms of the minimal generating set.
\end{Bemerkung}

The degree bound of the relations in this case is $12$. This allows us to
state the following conjecture.

\begin{Vermutung}
	$\nh^1(n,d) \le n(n+1)$.
\end{Vermutung}

\section{Proof}\label{proof}
The idea of the proof is to compute the Hilbert series of the candidate ideal.
For this one has to compute a Gr\"obner basis. In the initial setup there
are too many variables, so the Gr\"obner basis computation is too expensive.
But we can reduce the problem by the following homogeneous system of parameters.
\begin{Thm} [\cite{lopatin}] \label{hompar} The following elements form a 
	homogeneous system of parameters of $C_{33}$.\\
\begin{tabular}{cccc}
$\tr(X)$&$\tr(Y)$&$\tr(Z)$\\
$\tr(x^2)$&$\tr(xy)$&$\tr(xz)$\\ 
$\tr(y^2)$&$\tr(yz)$&$\tr(z^2)$\\ 
$\tr(x^3)$&$\tr(y^3)$&$\tr(z^3)$\\ 
$\tr(x^2y)-\tr(y^2z)-\tr(xz^2)$&$\tr(x^2z)-\tr(z^2y)$&$\tr(x^2z)-\tr(y^2x)$\\ 
&$\tr(xyz)-\tr(xzy)$&\\
$\tr(x^2y^2)$&$\tr(x^2z^2)$&$\tr(y^2z^2)$\\
\end{tabular}
\end{Thm}
The Theorem of Lopatin is much more general. Here we just picked one
special system of parameters of $C_{33}$.

For all elements of degree $\le 3$ is obvious, how they correspond to
the generators. For those of degree $4$ we get the relations
 \begin{eqnarray*}
  \tr(x^2y^2) &=& \frac{1}{6} t_1 t_2 + \frac{1}{3}t_4^2 + \frac{1}{3} t_{18},\\
  \tr(x^2z^2) &=& \frac{1}{6} t_1 t_3 + \frac{1}{3}t_5^2 + \frac{1}{3} t_{19},\\
  \tr(y^2z^2) &=& \frac{1}{6} t_2 t_3 + \frac{1}{3}t_6^2 + \frac{1}{3} t_{20}.\\
 \end{eqnarray*}
Because $t_1,\ldots,t_6$ are also elements of the homogeneous system of
parameters, we can replace the traces on the left by $t_{18}$,$t_{19}$ and
$t_{20}$ and get the following system of parameters:
\begin{eqnarray*}
	t_a,t_b,t_c &\\
	t_1,\ldots,t_9 &\\
	t_{10} - t_{13}-t_{14} &\\
	t_{11} - t_{15}&\\ 
	t_{11} - t_{12}&\\
	t_{17},\ldots, t_{20}.&\\
\end{eqnarray*}

If we divide out the homogeneous system of parameters we can eliminate some
variables by the following reductions.
\begin{eqnarray*}
	t_a,t_b,t_c,t_1,\ldots,t_9 &\rightsquigarrow& 0\\
	t_{10} &\rightsquigarrow& t_{13}+t_{14}\\
	t_{15} &\rightsquigarrow& t_{11}\\
	t_{12} &\rightsquigarrow& t_{11}\\
	t_{17},t_{18},t_{19},t_{20} &\rightsquigarrow& 0\\
\end{eqnarray*}

Let $J$ be the ideal given by \ref{main theorem} and $H$ be the homogeneous 
system of parameters above. Because the elements
of this homogeneous system of parameters are given in terms of traces, choose
the canonical preimage of these elements in $\C[T_a,\ldots,T_{45}]$ and denote
the ideal generated by $J$ and these elements by $J_H$. Then we get
the following inequality of the Hilbert series due to \cite{stanley}:
\begin{displaymath}
	\frac{H(\C[T_a,\ldots,T_{45}]/J_H,t)}{\prod (1-t^{d_i})} \ge 
	H(\C[T_a,\ldots,T_{45}]/J,t) \ge H(C_{33},t) = 
	\frac{H(C_{33}/H,t)}{\prod (1-t^{d_i})}.
\end{displaymath}
Here the products in the denominator are determined by the degrees of the
homogeneous parameter set.
The inequality in the middle holds, because we have the surjective map
$\pi \colon \C[T_a,\ldots,T_{45}]/J \rightarrow C_{33}$. Equality holds for our candidate
and $C_{33}$ if and only if the outer Hilbert series are equal. 
And these can be compared by comparing the nominators. 

The Hilbert series of $J_H$ can be computed by {\sc Singular} \cite{singular} since dividing out
the homogeneous system of parameters allows us to eliminate some variables.
Therefore {\sc Singular} can compute the Hilbert series of $\C[T_a,\ldots,T_{45}]/J_H$.
Since the Hilbert series of $C_{33}$ is known (see \cite{drensky:3x3def}) 
one gets that the elements given in \ref{main theorem} generate $\ker(\pi|_{C_E})$. 

Additional due to \cite{springer} every $\C$-basis $\C[T_a,\ldots,T_{45}]/J_H$
coming from homogeneous elements gives rise to a generating set of $C_{33}$
as a free $\C[H]$-module. From the Gr\"obner basis computation we get the following
corollary.

\begin{Korollar}
	$C_{33}$ is a free $\C[t_a,t_b,t_c,t_1,\ldots,t_9,t_{17},t_{18},t_{19},t_{20},t_{10}-t_{13}-t_{14},t_{11}-t_{15},t_{11}-t_{12}]$-module. The following elements and their divisors
form a basis of this module.
	\begin{longtable}{cccccccc}
$t_{11}^{2} t_{45}$&$t_{11} t_{13} t_{45}$&$t_{11} t_{14} t_{44}$&$t_{11} t_{14} t_{45}$&$t_{11} t_{16} t_{44}$&$t_{11} t_{16} t_{45}$&$t_{11} t_{22} t_{44}$&$t_{11} t_{23} t_{44}$\\
$t_{11} t_{23} t_{45}$&$t_{11} t_{24}$&$t_{11} t_{25} t_{45}$&$t_{11} t_{26} t_{44}$&$t_{11} t_{26} t_{45}$&$t_{11} t_{30}$&$t_{11} t_{31}$&$t_{11} t_{32} t_{45}$\\
$t_{11} t_{33}$&$t_{11} t_{34}$&$t_{11} t_{35} t_{45}$&$t_{11} t_{40}$&$t_{11} t_{41}$&$t_{11} t_{42}$&$t_{11} t_{43}$&$t_{13}^{2} t_{45}$\\
$t_{13} t_{14} t_{44}$&$t_{13} t_{14} t_{45}$&$t_{13} t_{16} t_{30}$&$t_{13} t_{16} t_{44}$&$t_{13} t_{16} t_{45}$&$t_{13} t_{22}$&$t_{13} t_{23} t_{35}$&$t_{13} t_{23} t_{45}$\\
$t_{13} t_{24}$&$t_{13} t_{25} t_{45}$&$t_{13} t_{26} t_{44}$&$t_{13} t_{26} t_{45}$&$t_{13} t_{31}$&$t_{13} t_{32} t_{45}$&$t_{13} t_{33}$&$t_{13} t_{34}$\\
$t_{13} t_{35} t_{45}$&$t_{13} t_{40}$&$t_{13} t_{42}$&$t_{14}^{2} t_{44}$&$t_{14}^{2} t_{45}$&$t_{14} t_{16} t_{31}$&$t_{14} t_{16} t_{32}$&$t_{14} t_{16} t_{43}$\\
$t_{14} t_{16} t_{44}$&$t_{14} t_{16} t_{45}$&$t_{14} t_{21}$&$t_{14} t_{22} t_{44}$&$t_{14} t_{22} t_{45}$&$t_{14} t_{23} t_{44}$&$t_{14} t_{23} t_{45}$&$t_{14} t_{24} t_{45}$\\
$t_{14} t_{25} t_{44}$&$t_{14} t_{25} t_{45}$&$t_{14} t_{26} t_{44}$&$t_{14} t_{26} t_{45}$&$t_{14} t_{30}$&$t_{14} t_{31} t_{45}$&$t_{14} t_{32} t_{45}$&$t_{14} t_{33}$\\
$t_{14} t_{34} t_{44}$&$t_{14} t_{34} t_{45}$&$t_{14} t_{35} t_{44}$&$t_{14} t_{35} t_{45}$&$t_{14} t_{41} t_{44}$&$t_{14} t_{41} t_{45}$&$t_{14} t_{42} t_{44}$&$t_{14} t_{42} t_{45}$\\
$t_{16}^{2} t_{43}$&$t_{16}^{2} t_{44}$&$t_{16}^{2} t_{45}$&$t_{16} t_{22} t_{44}$&$t_{16} t_{22} t_{45}$&$t_{16} t_{23} t_{44}$&$t_{16} t_{23} t_{45}$&$t_{16} t_{24}$\\
$t_{16} t_{25} t_{44}$&$t_{16} t_{25} t_{45}$&$t_{16} t_{26} t_{44}$&$t_{16} t_{26} t_{45}$&$t_{16} t_{27} t_{44}$&$t_{16} t_{27} t_{45}$&$t_{16} t_{28} t_{45}$&$t_{16} t_{29} t_{43}$\\
$t_{16} t_{29} t_{45}$&$t_{16} t_{30} t_{40}$&$t_{16} t_{31} t_{45}$&$t_{16} t_{32} t_{45}$&$t_{16} t_{33}$&$t_{16} t_{34} t_{44}$&$t_{16} t_{34} t_{45}$&$t_{16} t_{35} t_{45}$\\
$t_{16} t_{39}$&$t_{16} t_{40} t_{45}$&$t_{16} t_{41} t_{45}$&$t_{16} t_{42} t_{45}^{2}$&$t_{16} t_{44}^{2}$&$t_{16} t_{44} t_{45}$&$t_{21} t_{33} t_{45}$&$t_{21} t_{34} t_{45}$\\
$t_{21} t_{35} t_{45}$&$t_{21} t_{38}$&$t_{21} t_{40}$&$t_{21} t_{42} t_{45}$&$t_{21} t_{43}$&$t_{21} t_{44}$&$t_{21} t_{45}^{2}$&$t_{22} t_{33} t_{45}$\\
$t_{22} t_{34}$&$t_{22} t_{35} t_{45}$&$t_{22} t_{37}$&$t_{22} t_{39}$&$t_{22} t_{40}$&$t_{22} t_{41} t_{45}$&$t_{22} t_{42}$&$t_{22} t_{43}$\\
$t_{22} t_{45}^{2}$&$t_{23} t_{33}$&$t_{23} t_{34} t_{45}$&$t_{23} t_{35} t_{45}$&$t_{23} t_{38}$&$t_{23} t_{40}$&$t_{23} t_{41}$&$t_{23} t_{42} t_{45}$\\
$t_{23} t_{43}$&$t_{23} t_{45}^{2}$&$t_{24} t_{27} t_{45}$&$t_{24} t_{30}$&$t_{24} t_{31}$&$t_{24} t_{32}$&$t_{24} t_{40}$&$t_{24} t_{42}$\\
$t_{24} t_{43}$&$t_{24} t_{44}$&$t_{24} t_{45}^{2}$&$t_{25} t_{28} t_{45}$&$t_{25} t_{30}$&$t_{25} t_{32}$&$t_{25} t_{39}$&$t_{25} t_{40}$\\
$t_{25} t_{41}$&$t_{25} t_{42}$&$t_{25} t_{43}$&$t_{25} t_{45}^{2}$&$t_{26} t_{29} t_{45}$&$t_{26} t_{31}$&$t_{26} t_{32}$&$t_{26} t_{40}$\\
$t_{26} t_{41}$&$t_{26} t_{42}$&$t_{26} t_{43}$&$t_{26} t_{45}^{2}$&$t_{27} t_{40}$&$t_{27} t_{43}$&$t_{28} t_{41}$&$t_{28} t_{44}$\\
$t_{29} t_{42}$&$t_{29} t_{44}$&$t_{30} t_{41}$&$t_{30} t_{44}$&$t_{30} t_{45}^{2}$&$t_{31} t_{42}$&$t_{31} t_{44}$&$t_{31} t_{45}^{2}$\\
$t_{32} t_{42}$&$t_{32} t_{44}$&$t_{32} t_{45}^{2}$&$t_{33} t_{39}$&$t_{33} t_{40}$&$t_{33} t_{41}$&$t_{33} t_{44}$&$t_{33} t_{45}^{2}$\\
$t_{34} t_{40}$&$t_{34} t_{42}$&$t_{34} t_{43}$&$t_{34} t_{45}^{2}$&$t_{35} t_{41}$&$t_{35} t_{42}$&$t_{35} t_{43}$&$t_{35} t_{45}^{2}$\\
$t_{36} t_{44}$&$t_{37} t_{42}$&$t_{37} t_{43}$&$t_{37} t_{44}$&$t_{37} t_{45}$&$t_{38} t_{41}$&$t_{38} t_{44}$&$t_{38} t_{45}$\\
$t_{40} t_{42}$&$t_{40} t_{44}$&$t_{41} t_{42} t_{45}$&$t_{42}^{2} t_{45}$&$t_{43} t_{44}$&$t_{44}^{3}$&\\

	\end{longtable}
\end{Korollar}

\begin{Bemerkung}
	The first algorithm (section \ref{algorithms}) was implemented first in {\sc Maple} 
	\cite{maple} and later in {\sc Sage} \cite{sage}. That the relations given by this
	algorithm generate $\ker \pi|_{C_E}$ was confirmed using {\sc Singular} \cite{singular}
	by comparing the Hilbert series (section \ref{proof}).

	The second algorithm was also implemented in {\sc Sage}. 
\end{Bemerkung}

\bibliographystyle{alpha}
\newcommand{\etalchar}[1]{$^{#1}$}
\def\cdprime{$''$} \def\cprime{$'$}


\end{document}